\documentclass[12pt,draftcls,onecolumn]{IEEEtran}  

\IEEEoverridecommandlockouts                              

\usepackage{graphicx} 
\usepackage{epsfig} 
\usepackage{mathptmx} 
\usepackage{times} 
\usepackage{amsmath} 
\usepackage{amssymb}  
\usepackage{soul}
\usepackage[mathcal]{euscript} 
\usepackage{mathtools}
\usepackage{mathrsfs}

\usepackage{graphicx}
\graphicspath{{Figures/}} 
\usepackage{caption} 
\usepackage{subfig} 

\usepackage{url}

\usepackage{xcolor}
\definecolor{USred}{rgb}{0.74,0.1,0.1}
\definecolor{USblue}{rgb}{0.2,0.2,0.7}
\definecolor{green1}{cmyk}{0.82,0,1,0.3}

\definecolor{Royalblue}{cmyk}{1,0.30,0.2,0.2}
\definecolor{persiangreen}{rgb}{0.0, 0.65, 0.58}
\definecolor{ao}{rgb}{0.0, 0.5, 0.0}
\definecolor{amber(sae/ece)}{rgb}{1.0, 0.49, 0.0}

\newcommand{\RR}{\mathbb{R}}

\newcommand{\TT}{\mathbb{T}}
\newcommand{\ZZ}{\mathbb{Z}}
\newcommand{\EE}{\mathbb{E}}

\newcommand{\HH}{\mathbb{H}}

\newcommand{\D}{\mathcal{D}}

\newcommand{\V}{\mathcal{V}}
\newcommand{\G}{\mathcal{G}}
\newcommand{\I}{\mathcal{I}}
\newcommand{\J}{\mathcal{J}}

\newcommand{\N}{\mathcal{N}}
\newcommand{\Q}{\mathcal{Q}}

\newcommand{\argmin}{\operatornamewithlimits{argmin}}
\newcommand{\argmax}{\operatornamewithlimits{argmax}}
\renewcommand{\vec}{\boldsymbol}
\DeclareMathOperator{\tr}{tr}

\DeclareMathOperator{\diag}{diag}

\DeclareRobustCommand{\vect}[1]{
	\ifcat#1\relax
	\boldsymbol{#1}
	\else
	\mathbf{#1}
	\fi}

\DeclareRobustCommand{\rv}[1]{\vec{\mathsf{#1}}}

\newcommand\indep{\protect\mathpalette{\protect\independenT}{\perp}}
\def\independenT#1#2{\mathrel{\rlap{$#1#2$}\mkern2mu{#1#2}}}

\newtheorem{theorem}{Theorem}
\newtheorem{lemma}{Lemma}

\newtheorem{problem}{Problem}
\newtheorem{remark}{Remark}
\newtheorem{example}{Example}

\title{Data-driven Link Prediction over Graphical Models}

\author{Daniele Alpago, Mattia Zorzi, Augusto Ferrante 
	\thanks{}
	\thanks{D. Alpago, M. Zorzi and A. Ferrante are with the Department of Information Engineering, University of Padova, Padova, Italy; email:	 
		{\tt\small alpagodani@dei.unipd.it} (D. Alpago)
		{\tt\small zorzimat@dei.unipd.it} (M. Zorzi)
		{\tt\small augusto@dei.unipd.it} (A. Ferrante)}%
	\thanks{}%
}

\begin{document}
	
\maketitle


\begin{abstract}
The positive link prediction (PLP) problem is formulated in a system identification framework: we consider dynamic graphical models for auto-regressive moving-average (ARMA) Gaussian random processes. For the identification of the parameters, we model our network on two different time scales: a quicker one, over which we assume that the process representing the dynamics of the agents can be considered to be stationary, and a slower one in which the model parameters may vary. The latter accounts for the possible appearance  of new edges. The identification problem is cast into an optimization framework which can be seen as a generalization of the existing methods for the identification of ARMA graphical models. We  prove the existence and uniqueness of the solution of such an optimization problem and we propose a procedure to compute numerically this solution. Simulations testing the performances of our method are provided.
\end{abstract}



\section{Introduction}\label{sec:intro}
The widespread use of online social networks (e.g. Facebook, Twitter, YouTube, etc.) and the availability of a huge amount of high-quality data related to them, have raised an increasing interest on network analysis research. The nature of big data coming from social networks is rather complex due to the dynamical behavior of the network, responsible of changes in topological and non-topological features of the network over time. In this scenario, the so called \emph{link prediction problem} has been of great interest in the research community over the recent years \cite{Wang2015,liben2007, yang2015}. Given a certain network, link prediction problems can be divided in two categories: predicting edges that appear or disappear in the future and detecting edges in the current social network that are still unobserved. Link prediction techniques find important applications in recommending systems, in finding business collaborators in e-commerce networks \cite{Akcora2011}, academic social networks (describing co-authors relations) \cite{Wohlfarth2008}, security-related networks \cite{Huang2009}, e-mails networks \cite{Wang2007}, or gene-expression networks in biology \cite{Almansoori2012}. The most part of the research effort has been focused on the \emph{positive link prediction} (PLP) problem which consists either in detecting newly  appeared edges or in predicting unobserved links that are likely to appear, depending on the perspective we are adopting.

Most of the literature on PLP concerns the design of similarity measures between nodes that are unconnected according to the current network topology: an edge is added if the similarity measure between the corresponding nodes is sufficiently large with respect to some problem-dependent criterium. Many similarity measures have been proposed in the literature: common neighbors, Adamic/Adar, Katz, spreading activation, etc, most of which are designed either by taking into account the topological behavior exhibited in the past by the real network,  or on the basis of some properties the real network is expected to have, see \cite{Wang2015,liben2007} and references therein. For instance, the real network might posses the so-called ``small world'' property, meaning that most pairs of nodes are related through short chains and the network has an elevated clustering coefficient, a couple of unconnected nodes having several common neighbors should have a high similarity score. There are, however, some ``unfriendly prediction networks'', for which most of the similarity measures provide unsatisfying prediction accuracy, e.g. Facebook, UC Irvine, see \cite{gao2015link}.\\
In this paper we consider Gaussian graphical models whose  nodes are modeled as autoregressive, moving-average (ARMA) dynamical models and an edge between two nodes means that they are conditionally dependent given the other nodes of the network.  Our aim is to propose a new data-driven PLP paradigm leading to a similarity measure which exploits the current graphical model and some available noisy data. Hence, the properties on the underlying network are dictated from the data rather than from properties the network is expected to have. The development of such similarity measure is rephrased in the context of system identification with prior: we search the updated graphical model  which explains the data and that is as much close as possible to the current graphical model (i.e. the prior). The identification paradigm is then cast into an optimization problem for which we prove the existence and uniqueness of the solution. The similarity measure our method induces does not only exploits the topology of the current network (as the classic similarity measures do), but also the model of the current network making it more accurate for prediction.\\
It is worth noting that our contribution actually considers a \emph{detection} problem in which the prediction exploits both the past information (the current network) and a noisy piece of information (data) coming from the network where the new link has already appeared. The prediction scheme we propose is based on integrating the current network with  the data and is therefore ready for a recursive implementation. Starting from the time-series data, we also propose an online learning scheme for the parameters of a graphical model describing the network for which we want to detect new links. The modeling of our network has been carried on over two different time scales: a quicker one, over which we assume that the process representing the dynamics of the agents can be considered to be stationary, and a slower one in which the model parameters may vary.

Our work can be naturally cast in the context of identification of dynamic graphical models that are particularly interesting in modeling high-dimensional data \cite{Lauritzen1996,dahlhaus2000a,avventi2013,SongDhaVan2010,ZorzSep2016}. Particularly appropriate for big-data applications are in fact sparse graphical models that are characterized by a reduced number of model's parameters, thus reducing the risk of overfitting in the estimation scheme. The problem of identifying graphical models is a major topic in statistics and it has been tackled in various ways. Songsiri et al. \cite{SongVan2010}, for instance, proposed to carry out the identification of dynamic sparse graphical models associated to Gaussian autoregressive (AR) processes solving a regularized maximum likelihood problem. A Bayesian reformulation of the latter has been proposed in \cite{zorzi2019empirical}. Concerning the identification of \emph{ARMA graphical models}, we  mention \cite{avventi2013} in which the identification of Gaussian ARMA models has been addressed. Here, the MA part  is introduced as prior and it has a particular structure. The identification is carried over by employing a convex optimization approach solving a generalized moment problem like the ones studied in many papers over the recent years; see \cite{ByrnGusLind1998,ByrnGeoLind2000,GeoLind2003,georgiou2006,zorzi2014,ringh2016,zorzi2015interpretation,karlsson2013uncertainty}, just to name a few. We will show that our system identification paradigm with prior generalizes the ARMA graphical model identification framework of \cite{avventi2013}.

The paper is organized as follows: the notation and the background material on static and dynamic graphical models used throughout the paper is recalled in Section \ref{sec:notation}. Section \ref{sec:PLP} starts with the introduction of the PLP problem in terms of a regularized optimization problem. An interesting maximum likelihood interpretation and a recursive version of the optimization problem just mentioned are reported in Section \ref{sec:PLP} as well. In Section \ref{sec:existence} the existence and the uniqueness of the solution of the optimization problem previously set-up is proved. Section \ref{sec:sim} reports results of some numerical experiments meant to test the performances of the proposed method. Finally, in Section \ref{sec:conc} we draw the conclusions and analyze some possible further extensions of the present work.

\section{Notation and Background}\label{sec:notation}
\subsection{Notation}\label{subsec:notation}
Given a matrix $F$, $F^\top$ will denote the transpose of $F$, and $F^*$ its transpose-conjugate.
If $F$ is square of dimension $p$, $\tr(F),\,\det(F),\,F^{-1}$ stand for the trace, the determinant and the inverse of $F$, respectively; moreover, $\diag(F)$ denotes the $p$-dimensional vector whose entries are the diagonal elements of $F$ while $F\ge0$ and $F>0$ denote that $F$ is positive semi-definite or positive definite, respectively. $I_p$ denotes the $p\times p$ identity matrix.
We denote by $\HH_p$ the real vector space of Hermitian matrices of dimension $p\times p$ and by $L^\infty(\TT,\HH_p)$  the Banach space of essentially bounded functions defined on the unit circle $\TT:=\{e^{i\theta}:\,\theta\in[-\pi,\pi]\}$ and taking values in $\HH_p$.  For any function in $L^\infty(\TT,\HH_p)$, when it is clear from the context, we will drop the explicit dependence on $\theta$  and we  use the short-hand notation $\int\,F$ for
for the integral $\int_{-\pi}^\pi\,F(e^{i\theta}) \,\frac{d\theta}{2\pi}$.
In this paper we will deal in particular with functions in 
\begin{align*}
\mathcal{S}_p^+ :=\big\{ F\in L^\infty(\TT,\HH_p):\,F-\alpha\,I_p \ge0 \text{ a.e. on $\TT$},\, \exists\,\alpha>0\big\}
\end{align*}
namely, we will consider coercive, bounded functions defined on $\TT$. The operator $\mathsf{P}_\Omega:L^\infty(\TT,\HH_p)\to L^\infty(\TT,\HH_p)$ that maps a function $F$ to its projection $\mathsf{P}_\Omega(F)$ onto a certain support $\Omega\subseteq\{1,\dots,p\}\times\{1,\dots,p\}$, is defined as
\[
\left[\mathsf{P}_\Omega(F)\right]_{ij}=
\left\{
\begin{split}
&0\quad&\text{if}&\quad(i,j)\in\Omega^c\\
&F_{ij}\quad&\text{if}&\quad(i,j)\in\Omega
\end{split}
\right.
\]
where $\Omega^c$ denotes the complement of $\Omega$ in $\{1,\dots,p\}\times\{1,\dots,p\}$. We define  the following finite dimensional vector space of matrix pseudo-polynomials 
\[
\scalebox{0.9}{$\displaystyle
\mathcal{P}_{p,n}:=\left\{P\in L^\infty(\TT,\HH_p):\, P(e^{i\theta})=\sum_{k=-n}^n\,P_k\,e^{i\theta k},\, P_{-k}=P_k^\top\in\RR^{p\times p}\right\}
$}
\]
endowed with the norm
\begin{equation}\label{eq:normsp}
\|P\|_\mathcal{P}:=\int_{-\pi}^\pi\,|\nu_p(e^{i\theta})|\,\frac{d\theta}{2\pi},
\end{equation}
where $\nu_p(e^{i\theta}),\,\theta\in[-\pi,\pi]$, is the  eigenvalue of $P$ having maximum modulus. For a function $f$ defined in some metric space $X$ and taking values in $\RR$, we denote by $\text{epi}(f)$ its \emph{epigraph} namely, the subset of $X\times \RR$ defined by
\begin{equation*}
\text{epi}(f):=\left\{(x,\beta)\in X\times \RR:\,f(x)\le \beta\right\}.
\end{equation*}
The symbol $\EE[\cdot]$ denote the mathematical expectation.

\subsection{Preliminaries on Graphical models}\label{sec:graphmod}
A static undirected graphical model is a graph $\G=(V,E)$ associated to a certain $m$-dimensional Gaussian random vector $\rv{x}\sim\N(0,\Sigma)$, $\Sigma=\Sigma^\top>0$, having vertexes $V=\{1,\dots,m\}$ representing the components $\vec{x}_1,\dots,\vec{x}_m$ of $\rv{x}$ and edges $E\subset V\times V$ describing the conditional dependence relations between the components through the following equivalent relations
\begin{equation}\label{eq:cistatic}
(i,j)\notin E \quad\iff\quad \vec{x}_i \indep \vec{x}_j \,|\, \{\vec{x}_k\}_{k\ne i,j}\quad\iff\quad (\Sigma^{-1})_{ij}=0,
\end{equation}
see for instance \cite{Lauritzen1996} for further details. Following \cite{brillinger1996, dahlhaus2000a, avventi2013} it is possible to provide a dynamic counterpart of relation \eqref{eq:cistatic}. In fact, let $\rv{y}:=\{\rv{y}(t),\,t\in\ZZ\}$ be an $m$-dimensional, purely non-deterministic, Gaussian stationary process and let $\Phi(e^{i\theta})$ its power spectral density defined for $\theta\in[-\pi,\pi]$. 
For any index set $I\subset V$, define as $\mathcal{X}_I := \text{span}\{\vec{x}_j(t):\,j\in I,\,t\in\ZZ\}$ the closure of the set containing all the finite linear combinations of the variables $\vec{x}_j(t)$, $j\in I$, so that for any $i\ne j$, the notation
\begin{equation*}
\mathcal{X}_{\{i\}} \indep \mathcal{X}_{\{j\}} \,|\,\mathcal{X}_{V\setminus\{i,j\}}
\end{equation*}
means that for all $t_1,t_2$,  $\vec{x}_i(t_1)$ and $\vec{x}_j(t_2)$ are conditionally independent given the space linearly generated by $\{\vec{x}_k(t),\,k\in V\setminus\{i,j\}, t\in\ZZ\}$. One can prove that
\begin{equation}\label{eq:cidin1}
\mathcal{X}_{\{i\}} \indep \mathcal{X}_{\{j\}} \,|\,\mathcal{X}_{V\setminus\{i,j\}}\quad\iff\quad [\Phi(e^{i\theta})^{-1}]_{ij}=0,
\end{equation}
for any $\theta\in[-\pi,\pi]$. The latter generalizes the static relation \eqref{eq:cistatic}. Accordingly, we can construct the undirected graph $\G=(V,E)$ representing the conditional dependence relations between the components of the process $\rv{y}$ by defining the set of edges as follows:
\begin{equation}\label{eq:cidin2}
(i,j)\notin E \quad\iff\quad \mathcal{X}_{\{i\}} \indep \mathcal{X}_{\{j\}} \,|\,\mathcal{X}_{V\setminus\{i,j\}}.
\end{equation}
In this framework, the graph $\G$ is completely characterized by the inverse power spectral density of the process however, the information enclosed on the inverse spectrum $\Phi^{-1}$ is richer than the one carried by the couple $\G=(V,E)$.\\

\section{Parametric Positive Link Prediction}\label{sec:PLP}
In this section the positive link prediction (PLP) problem is introduced in the simplest possible setting. This will serve as a starting point to rephrase the PLP problem for dynamic systems, starting from undirected graphical models for Gaussian time-series recalled in the previous section. To the best of our knowledge, this point of view has not yet been considered in the link-prediction community; we think, however, that it could be an interesting way to look at the problem from a different angle. The main motivations come from the following considerations:
\begin{itemize}
	\item[-] it seems very reasonable that a similarity measure should take into account the dynamics of nodes and edges and the relation between the two (and thus not only the topology of the network) in predicting if a connection between them will appear;
	\item[-] we want to detect pair of nodes that are {\em directly} connected while avoiding to detect a link between nodes that exhibit correlated behaviors due to common connection with other nodes.
For this reason, our similarity score is not induced by dependence but rather by {\em conditional dependence}. As we have seen, this is equivalent to detect the vanishing entries in the inverse of the (matricial) spectral density. Notice that, by modifying a single entry of the spectral density the whole zero pattern of its inverse may change completely so that our similarity score accounts for the information distributed on the  whole networks and has therefore a {\em global nature}. This contrasts with other data driven approaches (see, e.g. \cite{Huang2009}) where the similarity score is induced by a {\em local} property
involving only the pair of nodes under scrutiny or, at most, a neighborhood of the target nodes. 
%
%
\end{itemize}
Dynamic conditional independence networks cover both aspects; in fact the presence of a link between two nodes is established on the basis of the \emph{whole} network's dynamics. The approach we propose is therefore quite different from the classic link prediction techniques: we consider a time-series model for the nodes dynamics and our similarity measure is computed not only from from the information on the current network topology but also from the current network model through the collected data. The following example aims to introduce the PLP problem.\\

\begin{example}\label{ex:examplePLP}
	The network in Figure \ref{fig:PLPex} represents five agents $a,\,b,\,c,\,d,\,e$, as nodes, and their reciprocal friendship relations at a certain time $\sigma>0$ as solid edges. Suppose that at time $\tau>\sigma$ agents $b$ and $d$ became friends with $c$, i.e. the dashed-red edges will appear in the network. The PLP problem consists in detecting these red edges based on the information on the network available at the previous time $\sigma$ and some noisy data at time $\tau$.
	\begin{figure}[h!]
		\centering
		\includegraphics[scale=1]{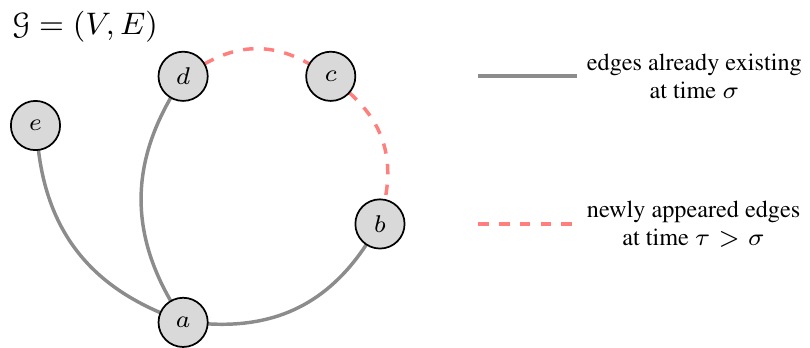}
		\caption{Network of Example \ref{ex:examplePLP}.}
		\label{fig:PLPex}
	\end{figure}
\end{example}

Hereafter we formalize our PLP problem. Let $\rv{y}:=\{\rv{y}(t),\,t\in\ZZ\}$ be an $m$-dimensional, Gaussian stochastic process and suppose that in a sufficiently small time interval $[\sigma-\alpha,\,\sigma+\alpha]$, for some $\sigma\in\ZZ_+$, $\alpha>0$, the process can be considered to be stationary \cite{dahlhaus2000b} so that it can be approximated by the parametric ARMA representation
\begin{equation}\label{eq:ARMA}
\sum_{k=0}^{n}\,A_j\,\rv{y}(t-k) = \sum_{k=0}^{n}\,B_j\,\rv{e}(t-k), \qquad \rv{e}(t)\sim\N(0,I_m),
\end{equation}
where $A_j,\,B_j\in\RR^{m\times m}$ for $j=0,1,\dots,m$, and $A_0=I_m$.
Denote with $R_k:=\EE[\rv{y}(t)\rv{y}(t-k)^\top]$, $k\in\ZZ$, its $k$-th covariance lag, so that the power spectral density of the process is the Fourier transform of the sequence $(R_k)_{k\in\ZZ}$,
\begin{equation}\label{eq:psd}
\Phi_\sigma(e^{i\theta}) = \sum_{k=-\infty}^{\infty}\,R_k\,e^{-i\theta k},\qquad R_{-k}=R_k^\top,\ \ \theta\in[-\pi,\pi].
\end{equation}
Accordingly, in the interval $[\sigma-\alpha,\,\sigma+\alpha]$ a graphical model $\G(\sigma) = (V,\,E_\sigma)$ associated to the process $\rv{y}$ is defined. In light of what we have said in the previous section, the edges of $\G(\sigma)$ are completely characterized by the support $\Omega_\sigma\subseteq\{(i,j):i,j=1,\dots,m\}$ of the inverse spectrum of the process $\Phi_\sigma^{-1}$, 
\[
\Omega_\sigma=\left\{(i,j):\, [\Phi_\sigma^{-1}(e^{i\theta})]_{ij}\ne0,\text{ for all }\theta\in[-\pi,\pi]\right\}.
\]
Notice that $\mathcal G(\sigma)$ represents the current network so both $\Phi_\sigma$ and $\Omega_\sigma$ are assumed to be known. Such model can be the result of previous studies on the network like for instance the so-called huge \emph{prospective investigations} carried out in medical research \cite{norat2005,gonzalez2006}, or it can be the outcome of a previous reliable estimation, more on this second case is developed in Section \ref{sec:recursive}. Let $\G(\tau)=(V,\,\Omega_\tau)$, $\tau>\sigma$, be the graphical model associated to the process $\rv{y}$ over another time interval $[\tau-\beta,\,\tau+\beta]$, $\beta>0$, where the process can be considered stationary and can be described by an ARMA representation as \eqref{eq:ARMA}, of course with different model parameters $\tilde{A}_j,\,\tilde{B}_j\in\RR^{m\times m}$ for $j=0,1,\dots,m$, and $\tilde{A}_0=I_m$. For our application it is very reasonable to assume that the model parameters do not change very much as the time interval change, see Remark \ref{rmk:slowchange} below. 

As above, $\G(\tau)$ corresponds to the support $\Omega_\tau$ of $\Phi_\tau^{-1}$: notice that both $\Phi_\tau^{-1}$ and its  support $\Omega_\tau$ are \emph{uknown}. Our aim is to estimate $\Phi_\tau$ corresponding to a certain support $\Omega_\tau$. 
Here $\G(\sigma)$ represents the known current network while $\G(\tau)$ plays the role of the new network, whose edges we want to predict. By virtue of relations \eqref{eq:cidin1} and \eqref{eq:cidin2}, the PLP scenario corresponds to the case in which $\Omega_\sigma\subseteq\Omega_\tau$, i.e. the network described by $\Phi_\tau^{-1}$ should have more edges than the network described by the prior spectrum $\Phi_\sigma^{-1}$.
Suppose now that $N$ observations $\mathsf{y}(1),\dots,\mathsf{y}(N)$ of the process over the interval $[\tau-\beta,\,\tau+\beta]$ are given as well. Natural interpretations of this \emph{nodes-related} information are statistics coming from questionnaires or interviews, very diffuse in the social networks framework. The available observations allow us to compute an estimate of the first $n+1$ covariance lags $R_0,\dots,R_n$ of $\rv{y}$ as
\begin{equation}\label{eq:estcov}
\hat{R}_k = \frac{1}{N}\sum_{t=k}^N\,\mathsf{y}(t)\,\mathsf{y}(t-k)^\top,\quad t\in[\tau-\beta,\,\tau+\beta],\quad k = 0,1,\dots,n.
\end{equation}
so that the truncated periodogram
\begin{equation}\label{eq:correl}
\hat{\Phi}_{n,\tau}(e^{i\theta}) = \sum_{k=-n}^n\,\hat{R}_k\,e^{-i\theta k},\qquad \hat{R}_{-k} = \hat{R}_k^\top,
\end{equation}
represents an estimate of the spectrum $\Phi_\tau$ based on the data. Notice that $\hat{\Phi}_{n,\tau}$ is not positive definite in general and, even if it were, $\hat{\Phi}_{n,\tau}^{-1}$ would not be sparse, in other words, it cannot be the estimate we are looking for. Our PLP problem can then be formalized as follows.\\

\begin{problem}\label{pb:spectralPLP}
	Given the prior power spectral density $\Phi_\sigma$ and the observations $\mathsf{y}(1),\dots,\mathsf{y}(N)$, estimate the spectrum $\Phi_\tau$ and the support $\Omega_\tau$ of its inverse, such that $\Omega_\tau\supseteq\Omega_\sigma$, $\Phi_\tau$ is as close  as possible to $\Phi_\sigma$, and the moments constraints $ \int e^{i\theta k} \,\Phi_\tau = \hat{R}_k$, are satisfied for $k=0,1,\dots,n$.\\   
\end{problem}

\begin{figure}[h!]
	\centering
	\includegraphics[scale=1]{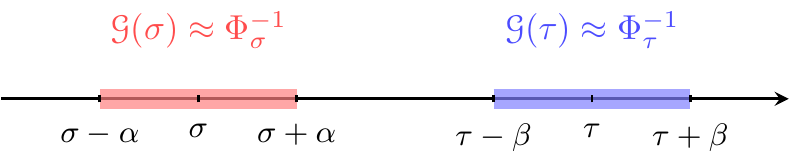}
	\caption{Pictorial clarification of the PLP set-up presented so far.}
	\label{fig:PLP}
\end{figure}

\begin{remark}\label{rmk:slowchange}
	Figure \ref{fig:PLP} attempts to give an intuitive explanation of the problem set-up. We stress the fact that, except for technical details, the framework is precisely the one informally presented in Example \ref{ex:examplePLP}. It is worth to note that the two intervals cannot be too far away from each other otherwise the model could be completely different depending on the interval that is selected. Considering the application this assumption is not as restrictive as it seems indeed, for instance, real life considerations tell us that it is very rare that edges on a social network change drastically over a relatively short period of time. This could happen after some ground-braking events that can be safely considered outliers (e.g. social media after the death of Michael Jackson, see for instance \cite{Kwak2010}).
\end{remark}

Having in mind the set-up just explained, we propose to reach the solution to Problem \ref{pb:spectralPLP} through the following steps:
\begin{enumerate}
	\item Consider a \emph{positive link selection problem}: a simplified version of the PLP problem in which the support $\Omega_\tau$ \emph{is known} and only the model dynamics described by $\Phi_\tau^{-1}$ needs to be estimated.
	\item Solve Problem \ref{pb:spectralPLP} by combining step 1) with the estimation of the support $\Omega_\tau$ from observations. 
\end{enumerate}
In order to ease the presentation, let us call $\Psi:=\Phi_\sigma$ the prior spectral density such that $\mathsf{P}_{\Omega_\sigma^c}(\Psi)=0$, $\Phi:=\Phi_\tau$ the spectrum we want to estimate and $\hat{\Phi}_n:=\hat{\Phi}_{n,\tau}$ the corresponding truncated periodogram \eqref{eq:correl}. Moreover, we assume that $\Phi,\Psi\in\mathcal{S}_m^+$ and that $\Psi$ is rational. It is worth noting that most of the results of this paper could be extended to a non-rational $\Psi$ at the price of technical complications that would obscure the presentation.

In light of the previous observations, we start by assuming the support $\Omega_\tau$ to be \emph{known}. Accordingly, Problem \ref{pb:spectralPLP} reduces to the estimation of the spectral density $\Phi$, whose inverse needs to have the specified support $\Omega_\tau$, given the prior $\Psi$ and the observations $\mathsf{y}(1),\dots,\mathsf{y}(N)$. This scenario is appropriate for setting-up the following link-selection problem:
\begin{equation}\label{eq:ISmin}
\begin{aligned}
\argmin_{\Phi\in\mathcal{S}_m^+} &\quad \D(\Phi\|\Psi)\\
\text{subject to } &\quad \mathsf{P}_{\Omega_\tau}\left(\int e^{i\theta k} \,\Phi - \hat{R}_k\right)=0,\qquad k=0,1,\dots,n.
\end{aligned}
\end{equation} 
where
\begin{equation}\label{eq:ISdist}
\D(\Phi\|\Psi):=\frac{1}{4\pi}\int_{-\pi}^{\pi}\,-\log\det\Phi+\log\det\Psi+\tr(\Psi^{-1}\Phi)\,d\theta-m
\end{equation}
is the multivariate form of Itakura-Saito pseudo-distance among power spectral densities $\Phi$ and $\Psi$, employed in signal processing \cite{basseville1989}. Computational-complexity arguments and rationality constraints on the spectra motivate the choice of \eqref{eq:ISdist} as cost functional to be minimized in these kind of problems, see \cite{FerMasPav2012,zorzi2014multivariate} for a deeper discussion on that.\\

\begin{remark}
	It is worth noticing that the support's constraints on the moments in Problem \eqref{eq:ISmin} are the natural way to make the problem accounting the information provided by the observations $\mathsf{y}(1),\dots,\mathsf{y}(N)$, that has to be used in the identification procedure, together with the prior $\Psi$. In fact, we do not impose that $\Phi$ matches entirely the estimated moments $\hat{R}_k$, $k=0,1,\dots,n$. The condition in (\ref{eq:ISmin}) has indeed the following interpretation: we believe that only the covariances between the nodes in $\Omega_\tau$ are reliable.
\end{remark}
	
Starting from \cite{dempster1972, burg1975}, Dempster problems have been developed, progressively generalized and adapted to different frameworks in a blooming stream of literature including \cite{ByrnGusLind1998,ByrnGeoLind2000,GeoLind2003,georgiou2006,zorzi2014,ringh2016,ByrEnqLin2001,KarLinRin2016,KarGeoLin2010,Geo2005}. However, generalized Dempster problems in graphical-models scenarios always involve maximum entropy problems in which some kind of entropic functional is maximized under linear constraints \cite{avventi2013,SongDhaVan2010,SongVan2010,ZorzSep2016,AlpZorzFer2018}. In this respect, Problem \eqref{eq:ISmin} in which a \emph{relative entropy functional is minimized}, appears to
be an original approach for graphical models applications.\\
As in Dempster's setting \cite{dempster1972,PavFer2013}, the natural representation of our problem is given in terms of inverse power spectral densities. Duality theory proves to be the right tool for re-parametrize the problem. The Lagrangian for Problem \eqref{eq:ISmin} is
\begin{equation*}
\mathcal{L}(\Phi,\tilde{\Lambda}) = \D(\Phi\|\Psi) - \sum_{k=0}^{n}\left<\mathsf{P}_{\Omega_\tau}\left(\hat{R}_k-\frac{1}{2\pi}\int\nolimits_{-\pi}^{\pi}\Phi\,e^{i\theta k}\,d\theta\right),\,\tilde{\Lambda}_k\right>,
\end{equation*}
where $\tilde{\Lambda}=[\tilde{\Lambda}_0\quad\tilde{\Lambda}_1\quad\cdots\quad\tilde{\Lambda}_n]$, $\tilde{\Lambda}_{-k}=\tilde{\Lambda}_k^\top\in\RR^{m\times m}$, encloses the Lagrange multipliers. Recalling that the projection operator is self-adjoint, we can rewrite the Lagrangian as
\begin{equation*}
\mathcal{L}(\Phi,\Lambda) = \D(\Phi\|\Psi) - \sum_{k=0}^{n}\left<\hat{R}_k-\frac{1}{2\pi}\int\nolimits_{-\pi}^{\pi}\Phi\,e^{i\theta k}\,d\theta,\,\Lambda_k\right>,
\end{equation*}
where $\Lambda=[\Lambda_0\quad\Lambda_1\quad\cdots\quad\Lambda_n]$ with $\Lambda_k:=\mathsf{P}_{\Omega_\tau}(\tilde{\Lambda}_k),\,k=0,1,\dots,n$, are the new multipliers. Moreover, from 
\[
\hat{R}_k = \frac{1}{2\pi}\int\nolimits_{-\pi}^\pi\hat{\Phi}_n\,e^{i\theta k}d\theta,\qquad k=0,1,\dots,n,
\]
we get
\begin{align*}
\mathcal{L}(\Phi,\Lambda) &= \D(\Phi\|\Psi) -\frac{1}{2\pi}\int_{-\pi}^\pi\sum_{k=0}^{n}\tr\left[\Lambda_k^\top\left(\hat{\Phi}_n-\Phi\right)e^{i\theta k}\right]\,d\theta\\
&=\frac{1}{4\pi}\int_{-\pi}^\pi-\log\det\Phi+\log\det\Psi+\tr\left(\Psi^{-1}\Phi\right)\\
&+2\sum_{k=0}^n\tr\left(\Lambda_k^\top\Phi\right)e^{i\theta k}-2\sum_{k=0}^n\tr\left(\Lambda_k^\top\hat{\Phi}_n\right)\,d\theta - m.
\end{align*}
By defining the pseudo-polynomial
\[
Q(e^{i\theta}) := \sum_{k=-n}^n\,Q_k\,e^{-i\theta k},
\ \ \text{with}\ \  Q_k:=
\left\{
\begin{split}
2\,\Lambda_0\quad\text{if}\quad k=0,\\
\Lambda_k^\top \quad\text{if}\quad k\ne0,
\end{split}   
\right.
\]
we can re-parametrize the Lagrangian as
\begin{align*}
\mathcal{L}(\Phi,Q) =
\frac{1}{4\pi}\int_{-\pi}^\pi&-\log\det\Phi+\log\det\Psi\\
&+\tr\left[\left(\Psi^{-1}+Q\right)\Phi\right] - \tr\left(Q\,\hat{\Phi}_n\right)\,d\theta - m	
\end{align*}
Notice that by construction, $Q=Q(e^{i\theta})$ has support $\Omega_\tau$, i.e. $\mathsf{P}_{\Omega_\tau^c}(Q)=0$. Thanks to the monotonicity of the integral, it suffices to minimize the integrand (strictly convex in $\Phi$),
\begin{align*}
f(\Phi,Q):=&-\log\det\Phi+\log\det\Psi\\
&+\tr\left[\left(\Psi^{-1}+Q\right)\Phi\right] - \tr\left(Q\,\hat{\Phi}_n\right).
\end{align*}
If $\delta f(\Phi,Q;\delta\Phi)$ denotes the G\^ateaux derivative of $f$ in a certain direction $\delta\Phi\in L^\infty(\TT,\HH_m)$ we have that
\begin{equation*}
\delta f(\Phi,Q;\delta\Phi)=0,\qquad\forall\,\delta\Phi\in L^\infty(\TT,\HH_m),
\end{equation*}
if and only if
\begin{equation}\label{eq:optcond}
	\Phi = \Phi_o = \left(\Psi^{-1}+Q\right)^{-1},
\end{equation}
provided that $Q\in\mathcal{P}_{m,n}$ is chosen such that $\Psi^{-1}+Q>0$ on $[-\pi,\pi]$. Observe that $\Phi_o^{-1}$ computed as in \eqref{eq:optcond}, has support $\Omega_\sigma\cup\Omega_\tau=\Omega_\tau$ which agrees with the positive link prediction set-up. Accordingly, 
the dual problem readily follow
\begin{equation}\label{eq:ISminDual}
\begin{aligned}
\argmax_{Q\in\mathcal{Q}_\Psi^+} &\quad \frac{1}{4\pi}\int_{-\pi}^\pi\log\det\left(\Psi^{-1}+Q\right)-\tr\left(Q\,\hat{\Phi}_n\right)\,d\theta\\
\text{subject to } &\quad \mathsf{P}_{\Omega_\tau^c}\left(Q\right)=0,
\end{aligned}
\end{equation}
where
\[
\mathcal{Q}_\Psi^+ :=\left\{Q\in\mathcal{P}_{m,n}:\,\Psi^{-1}(e^{i\theta})+Q(e^{i\theta})>0,\,\forall\,\theta\in[-\pi,\pi]\right\}
\]
is the domain of optimization, open and unbounded. In view of (\ref{eq:optcond}) and (\ref{eq:ISminDual}) the solution to Problem (\ref{eq:ISmin}) is such that its inverse has support $\Omega_\tau$, i.e. we have more (conditional) dependencies between the variables in the new network. We refer the reader to Section \ref{sec:existence}, Theorem \ref{thm:exuqDemp}, for the discussion of existence and uniqueness of the solution of Problem \eqref{eq:ISmin}.\\
	
\begin{remark} \label{remark-avventi}
It is interesting to note that Problem (\ref{eq:ISmin}) solves the problem of identifying an ARMA graphical model with topology $\Omega_\tau$ provided that the prior spectral density $\Psi$ corresponds to a graphical model having topology $\Omega_\sigma\subseteq \Omega_\tau$. It is worth noting that a similar problem has been addressed in \cite{avventi2013} with the following formulation:
\begin{equation}\label{pb_avvventi}
\begin{aligned}
\argmax_{\Phi\in\mathcal{S}_m^+} &\quad -\tr\int \psi\,\log( \psi\, \Phi^{-1})\\
\text{subject to } &\quad \mathsf{P}_{\Omega_\tau}\left(\int e^{i\theta k} \,\Phi - \hat{R}_k\right)=0,\qquad k=0,1,\dots,n,
\end{aligned}
\end{equation}
where $\psi$ is a scalar spectral density a priori known. Notice that (\ref{pb_avvventi}) is equivalent to minimize
\[
\D_{KL}(\psi\,I_m\| \Phi ):=\tr\left\{\int \psi \log( \psi\,\Phi^{-1}) -\psi I_m+\Phi\right\}
\] 
where the last two terms are constant, as $\psi$ is known and the moments constraint in (\ref{pb_avvventi}) ensures that $\tr \int \Phi=\tr(\hat R_0)$. Using arguments similar to the ones in \cite{zorzi2014rational}, it is not difficult to see that $\D_{KL}(\psi I_m\| \Phi )$ is actually a pseudo-distance between $\psi I_m$ and $\Phi$ ($\D_{KL}(\psi I_m\| \Phi )\geq 0$ with equality if and only if $\psi I_m=\Phi$), representing the natural extension of the {\em Kullback-Leibler} divergence between multivariate power spectral densities in which the first argument has the particular structure $\psi I_m$. Therefore, our problem is in the same spirit of Problem (\ref{pb_avvventi}). In particular, if we take $\Psi=I_m$ in (\ref{eq:ISmin}) and $\psi=1$ in (\ref{pb_avvventi}), then the two problems do coincide, i.e. they maximize the entropy rate of the process. In our setting, however, the scalar prior $\psi I_m$ would correspond to a graphical model with disconnected nodes. In principle, one could extend Problem (\ref{pb_avvventi}) to a prior corresponding to a graphical model having topology $\Omega_\sigma\subseteq \Omega_\tau$. However, in this case, the variational analysis cannot not be carried out.\\
\end{remark}

The second step of our solution consists in combining Problem \eqref{eq:ISminDual} with the estimation of the support $\Omega_\tau$ in order to obtain an optimization problem for the solution of Problem \ref{pb:spectralPLP}. Following \cite{SongVan2010,ZorzSep2016}, we propose to perform this step by resorting to a regularized version of Problem \eqref{eq:ISminDual}, namely\footnote{Notice that in order to deal with convex functions, 
instead of maximizing the  objective function (as in  Problem \eqref{eq:ISminDual}) we are now minimizing the opposite function (multiplied by a factor 2).}
\begin{equation}\label{eq:regISminDual}
\argmin_{Q\in\mathcal{Q}_\Psi^+} \quad \J_\Psi(Q)
\end{equation}   
where
\begin{align*}
\J_\Psi(Q) :&= \int\,\left[-\log\det\left(\Psi^{-1}+Q\right)+\tr\left(Q\,\hat{\Phi}_n\right)\right] +\lambda\,h_\sigma^\infty(Q)\\
&=\int\,\tr\left[Q\,\hat{\Phi}_n-\log\left(\Psi^{-1}+Q\right)\right] +\lambda\,h_\sigma^\infty(Q)
\end{align*}
and 
\[
h_\sigma^\infty(Q) = \sum_{(h,k)\in\I_\sigma}\max\left\{ |(Q_0)_{hk}|,\max_{j=1,\dots,n}|(Q_j)_{hk}|,\max_{j=1,\dots,n}|(Q_j)_{kh}|\right\},
\]
with $\I_\sigma := \left\{(h,k)\in V\times V\setminus\Omega_\sigma:\,k>h\right\}$, plays the role of the $\ell^1$-norm used to induce sparsity on vectors and it has been proposed in \cite{SongVan2010} for inducing group-sparsity to $Q_0,\dots,Q_n$. It is worth noticing that the sparsity-inducing regularization, tuned by the parameter $\lambda>0$, acts \emph{only} on the elements of the $Q_k$s in positions that are not contained in the support $\Omega_\sigma$, according to the fact that $\Omega_\sigma\subseteq\Omega_\tau$. This allow us to reduce the bias introduced by the regularization in the estimation procedure. Indeed, regularization here is used to decide whether an edge is present or not. Since we already know the presence of the edges in $\I_\sigma^c$, we do not need regularization for them. The proof of the existence and the uniqueness of the solution to Problem \eqref{eq:regISminDual} is devoted to Section \ref{sec:existence}.

We are aware that at first sight, our set-up may seems far apart with respect to the typical settings adopted in the link prediction community. However, if we have a closer look to the problem we readily find out that we are actually just proposing a different choice of \emph{score matrix} to decide how much a pair of nodes are inclined to get connected. Our approach suggests that a suitable \emph{similarity measure} should be an indicator of the conditional dependence between the variables. In order to introduce such a measure, we define the \emph{partial coherence} of the predicted spectrum $\Phi_\tau^{-1}$ as
	\begin{equation}\label{eq:partch}
	\Gamma_\tau(e^{i\theta}) := \diag[\Phi_\tau(e^{i\theta})]^{1/2}\,\Phi_\tau^{-1}(e^{i\theta})\,\diag[\Phi_\tau(e^{i\theta})]^{1/2},
	\end{equation}
	for all $\theta\in[-\pi,\pi]$. This is a standard tool the frequency-domain analysis of time series and it measures the dependence between two time series after removing the linear time invariant effects of a the other series \cite{brillinger1996,dahlhaus2000a,SongVan2010}. The similarity measure that is naturally induced by our approach in order to rate an edge is therefore 
	\begin{equation*}\label{eq:simmes}
	(G_\tau)_{ij} := \sqrt{\int_{-\pi}^\pi\,|[\Gamma_\tau(e^{i\theta})]_{ij}|^2\,d\theta}, 
    \quad (i,j)\in V\times V\setminus\Omega_\sigma,
	\end{equation*}
	and the matrix $G_\tau:=[(G_\tau)_{ij}]$ represents our score matrix. As the score matrix $G_\tau$ will have some small entries but, in general, will not be exactly sparse, a thresholding procedure is needed in order to obtain an estimate of the support $\Omega_\tau$ that defines a network topology according to relations \eqref{eq:cidin1}-\eqref{eq:cidin2}. More precisely, we will consider the edge $(i,j)\in V\times V\setminus\Omega_\sigma$ to be in $\Omega_\tau$ only if its score is greater than a certain threshold $t_r>0$ (to be suitably selected) namely,
	\[
         (i,j)\in\Omega_\tau\quad\iff\quad (G_\tau)_{ij}>t_r.
    \]
    To conclude, it is interesting to note that given a pair $(i,j)\in V\times V\setminus\Omega_\sigma$, one can get a straightforward interpretation of $(G_\tau)_{ij}$ in term of best (linear) predictors of $\rv{y}_i$ and $\rv{y}_j$, namely    
	\[
	   \hat{\rv{y}}_i(t) := \EE\left[\rv{y}_i\,\big\vert\,\mathcal{Y}_{V\setminus\{i,j\}}\right]\quad\text{and}\quad \hat{\rv{y}}_j(t) := \EE\left[\rv{y}_j\,\big\vert\,\mathcal{Y}_{V\setminus\{i,j\}}\right]
	\]
	respectively, where $\mathcal{Y}_{V\setminus\{i,j\}}:=\text{span}\{\rv{y}_k(t):\,k\in V\setminus\{i,j\}\}$. Then $(G_\tau)_{ij}$ is related to the correlation between the estimation errors $\rv{\epsilon}_i(t):=\rv{y}_i(t)-\hat{\rv{y}}_i(t)$ and $\rv{\epsilon}_j(t):=\rv{y}_j(t)-\hat{\rv{y}}_j(t)$ by noticing that
	\[
	 (G_\tau)_{ij}=
	 \left\lVert\frac{\Phi_{\rv{\epsilon}_i,\rv{\epsilon}_j}}{\sqrt{\Phi_{\rv{\epsilon}_i}\Phi_{\rv{\epsilon}_j}}}\right\lVert_2,
	\]
	where $\Phi_{\rv{\epsilon}_i}$, $\Phi_{\rv{\epsilon}_j}$ are the spectra of $\rv{\epsilon}_i$, $\rv{\epsilon}_j$ and 
$\Phi_{\rv{\epsilon}_i,\rv{\epsilon}_j}$ is the corresponding cross spectrum.

\subsection{Maximum Likelihood Interpretation}\label{sec:maxlik}
In this section we show that Problem \eqref{eq:ISminDual} has a nice interpretation as (regularized) maximum likelihood problem. This interpretation is based on an frequency approximation of the likelihood function of a Gaussian sample, firstly introduced by Whittle \cite{whittle1953,whittle1954} for scalar stationary processes. Over the years, the so-called Whittle likelihood approximations have been generalized to multivariate stationary processes \cite{DunsHann1976} and also extended to the non-stationary case \cite{dahlhaus2000b}.\\
Consider process \eqref{eq:ARMA} having power spectral density $\Phi$ as in \eqref{eq:psd}, and let $\hat{\Phi}_n$ as in \eqref{eq:correl} be its \emph{truncated} periodogram computed on the basis of $N$ given observations of the process $\mathsf{y}(1),\dots,\mathsf{y}(N)$. The assumptions on model \eqref{eq:ARMA} that guarantee the following to hold are to some extent classical, however, they require a quite technical presentation that the interested reader can find in \cite[Sec. 2]{DunsHann1976}. Let $\vect{y}_N=[\mathsf{y}(1)^\top\,\cdots\,\mathsf{y}(N)^\top]^\top$, then the $N$-dimensional probability density $p(\vect{y}_N;\Phi):= p(\mathsf{y}(1),\dots,\mathsf{y}(N);\Phi)$ of the random variables $\mathsf{y}(1),\dots,\mathsf{y}(N)$ has the well-known form
\[
p(\vect{y}_N;\Phi)=\dfrac{1}{\sqrt{(2\pi)^N\det\vect{T}_N(\Phi)}}\exp\left\{-\frac{1}{2}\vect{y}_N^\top\,\vect{T}_N(\Phi)^{-1}\vect{y}_N\right\}
\]
where $\vect{T}_N(\Phi)=\EE[\vect{y}_N\,\vect{y}_N^\top]$ is the $mN\times mN$ Toeplitz matrix whose $(h,k)$-th block is defined as 
\[
T(\Phi)_{hk}  = \int_{-\pi}^\pi\,\Phi(e^{i\theta})e^{i(h-k)\theta}\,\frac{d\theta}{2\pi} =R_{h-k},\,\quad 1\le h,k\le N+1.
\]
The corresponding negative log-likelihood (up to scaling factors and constant terms) is
\begin{equation}\label{eq:neglik}
\tilde{\ell}_N(\Phi)=\frac{1}{N}\log\det\vect{T}_N(\Phi) + \frac{1}{N}\vect{y}_N^\top\vect{T}_N(\Phi)^{-1}\vect{y}_N.
\end{equation}
Various frequency approximations of \eqref{eq:neglik} may be introduced \cite{DunsHann1976,dzhaparidze1986}. In order for us to define one of those, we introduce the discrete Fourier transform of the data
\begin{equation*}
Y_N(e^{i\theta})=\frac{1}{\sqrt{N}}\sum_{p=1}^N \mathsf{y}(p)\,e^{-i\theta p}
\end{equation*}
such that
\begin{equation*}
\hat{\Phi}_N(e^{i\theta}) = Y_N(e^{i\theta})Y_N(e^{i\theta})^*
= \sum_{k=-(N-1)}^{N-1}\,\hat{R}_k\,e^{-i\theta k},
\end{equation*}
where $\hat{R}_k$, such that $\hat{R}_{-k}=\hat{R}_k^\top$, $k=0,1,\dots,N-1$,  computed as in \eqref{eq:estcov}, is the periodogram of the process $\rv{y}$. With these definitions it can be shown \cite{DunsHann1976} that $\left(\tilde{\ell}_N-\ell_\infty\right)\to0$ almost surely as $N\to\infty$, where
\begin{equation}\label{eq:neglikasy}
\ell_\infty(\Phi)=\int_{-\pi}^{\pi}\,-\log\det\Phi(e^{i\theta})^{-1} + \tr\left[\Phi(e^{i\theta})^{-1}\,\hat{\Phi}_N(e^{i\theta})\right]\,\frac{d\theta}{2\pi}.
\end{equation}
In this sense Problem \eqref{eq:ISminDual} can be interpreted, at least asymptotically, as a (regularized) maximum likelihood problem, in which we have to minimize $\ell_\infty(\Phi)$ with $\Phi$ belonging to the parametric family
\begin{equation*}
\mathfrak{P}_{\Psi,Q}:=\left\{\Phi=\left(\Psi^{-1}+Q\right)^{-1}:\,Q\in\Q_\Psi^+ \right\},
\end{equation*}
for a given $\Psi\in\mathcal{S}_m^+$. Notice that
	\[
	Q(e^{i\theta})\,\hat{\Phi}_N(e^{i\theta}) = \sum_{k,h=-(N-1)}^{N-1}\,Q_k\,\hat{R}_h\,e^{-i\theta(k+h)}, 
	\]
	where $Q_k=0$ for $k=n+1,\dots,N-1$, because $Q\in\mathcal{P}_{m,n}$. Hence,
	\[
	Q(e^{i\theta})\,\hat{\Phi}_N(e^{i\theta}) = \sum_{k,h=-n}^{n}\,Q_k\,\hat{R}_h\,e^{-i\theta(k+h)}=Q(e^{i\theta})\,\hat{\Phi}_n(e^{i\theta})
	\]
	where $\hat{\Phi}_n$ is the \emph{truncated} periodogram \eqref{eq:correl} of the process. Accordingly, minimizing \eqref{eq:neglikasy} over $\mathfrak{P}_{\Psi,Q}$ is equivalent to minimize 
\begin{equation}\label{eq:likefin}
\int_{-\pi}^{\pi}\,-\log\det\left(\Psi^{-1}+Q\right) + \tr\left(Q\,\hat{\Phi}_n\right)\,\frac{d\theta}{2\pi}
\end{equation}
over $\Q_\Psi^+$. Save for the regularization term, \eqref{eq:likefin} is precisely the index $\J_\Psi$ of \eqref{eq:ISminDual}.\\

\subsection{Recursive Positive Link Prediction}\label{sec:recursive}
\begin{figure}[h!]
	\centering
	\includegraphics[scale=1.05]{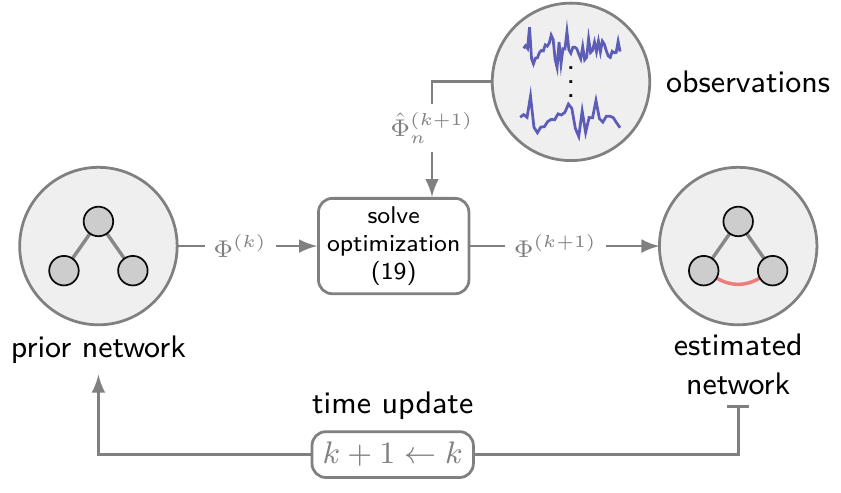}
	\caption{Iteration $k$ of recursive DPLP scheme.}
	\label{fig:iter}
\end{figure}

Our approach can be easily embedded in a recursive approach that perfects the previous estimation given new available data. Suppose that we want to study a certain network over the time interval $[0,\,T]$. We refer to the time window $[k-\alpha,k+\alpha]$, $\alpha>0$, $k>0$, as the time window $\tau=k$, see the beginning of this section. Moreover, we assume that at time $\tau=0$ a prior information concerning the whole network is summarized in the spectrum $\Phi_0$ whose inverse has support $\Omega_0$ and observations $\mathsf{y}^{(1)}(1),\dots,\mathsf{y}^{(1)}(N_1)$ related to the present time window $\tau=1$ are available. Then we can find $Q_1$ by solving Problem \eqref{eq:regISminDual} so that the predicted network at time $\tau=1$ is given by $\Phi_1^{-1}=\Phi_0^{-1}+Q_1$. The same reasoning applied at time $k+1$ leads to
\begin{equation}\label{eq:recopt}
	Q_{k+1}=\argmin_{Q\in\mathcal{Q}_{\Phi_k}^+} \quad \J_{\Phi_k}(Q)
\end{equation}
from which 
\[
  \Phi_{k+1}^{-1}=\Phi_k^{-1}+Q_{k+1} = \Phi_0^{-1}+\sum_{l=1}^{k+1}Q_l(e^{i\theta}),
\]
$Q_l\in \mathcal{P}_{m,n}$, $l=1,\dots,k+1$, being $\Phi_k$ the output for the estimation procedure at step $k$. Therefore, Problem (\ref{eq:recopt}) is equivalent to 
\begin{equation*}
	Q_{k+1}=\argmin_{Q\in\mathcal{Q}_{\Phi_0}^+} \quad \tilde \J_{\Phi_0,k}(Q),
\end{equation*} 
where
\[
 \tilde \J_{\Phi_0,k}(Q) := \frac{1}{4\pi}\int_{-\pi}^\pi\log\det\left(\Phi_0^{-1}+Q\right)-\tr\left(Q\,\hat{\Phi}_{n,k}\right)\,d\theta.
\]
What we just saw implies that  the upper-bound on \emph{Mcmillan degree} of $\Phi_k$ is constant and finite with respect to $k$, more precisely 
\[
 \deg(\Phi_k)\leq\deg(\Phi_0)+n, \; \; \forall\, k>0.
\]
It is worth noticing that, thanks to the latter upper-bound, the complexity of the model is guaranteed not to explode even if the number of iterations gets considerably high (potentially infinite). The corresponding iterative scheme is represented in Figure \ref{fig:iter}:
\begin{itemize}
	\item[-] \emph{Initialization}. Initial prior $\Phi_0$ and observations $\mathsf{y}^{(1)}(1),\dots,\mathsf{y}^{(1)}(N_1)$.
	\item[-] \emph{Iteration}. Given $\mathsf{y}^{(k+1)}(1),\dots,\mathsf{y}^{(k+1)}(N_{k+1})$ for $k=1,2,\dots,T$ compute
$Q_{k+1}$ by solving (\ref{eq:recopt}).
	
\end{itemize}


\section{Existence and Uniqueness of the Solution}\label{sec:existence}
This section is devoted to the proof of the existence and uniqueness of the solution for the regularized dual Problem \eqref{eq:regISminDual}. The idea of the proof is the following. We first extend the definition of $\J_\Psi(Q)$ to the boundary of $\Q_\Psi^+$ and show that this extended function do admit a unique minimum there. We then show that this minimum cannot be on the boundary.
Let
\begin{align*}
\partial\Q_\Psi^+ := \bigg\{Q\in\mathcal{P}_{m,n}:\,&\Psi^{-1}(e^{i\theta})+Q(e^{i\theta})\ge0 \text{ and singular},\\ &\exists\,\theta\in[-\pi,\pi]\bigg\}
\end{align*}
be the boundary of $\Q_\Psi^+$ and introduce the sequence of functions $\left(\J_\Psi^n\right)_{n\ge1}$ over 
\begin{align*}
\bar{\Q}_\Psi^+ :&= \Q_\Psi^+ \cup \partial\Q_\Psi^+\\
&=\left\{Q\in\mathcal{P}_{m,n}:\,\Psi^{-1}(e^{i\theta})+Q(e^{i\theta})\ge0,\,\exists\,\theta\in[-\pi,\pi]\right\}
\end{align*}    
defined as
\begin{align*}
\J^n_\Psi(Q) :=&\int_{-\pi}^\pi\tr\left[Q\,\hat{\Phi}_n-\log\left(\Psi^{-1}+Q+\frac{1}{n}I_m\right)\right]\,\frac{d\theta}{2\pi}\\&+\lambda\,h_\sigma^\infty(Q)
\end{align*}  
for $n\ge1$. Notice that $\J_\Psi^n$ is a strictly-convex and continuous function on $\bar{\Q}_\Psi^+$, therefore $\text{epi}\left(\J_\Psi^n\right)$ is convex and closed on $\bar{\Q}_\Psi^+\times\RR$ for any $n\ge1$. Moreover, the sequence $\left(\J_\Psi^n\right)_{n\ge1}$ is monotonically increasing so the pointwise limit 
\[
\J_\Psi^\infty(Q):=\lim_{n\to\infty}\,\J_\Psi^n(Q),
\]
exists and it is a strictly-convex and continuous function over $\bar{\Q}_\Psi^+$ because $\text{epi}\left(\J_\Psi^\infty\right)=\cap_{n\ge1}\text{epi}\left(\J_\Psi^n\right)$. 

First of all we ensure that $\J_\Psi^\infty$ coincides with the original $\J_\Psi$ on its domain of definition $\Q_\Psi^+$:
\begin{lemma}
The following relation holds on $\Q_\Psi^+$:
\begin{equation}
\J_\Psi^\infty\equiv \J_\Psi.
\end{equation}
\end{lemma}
\begin{IEEEproof}
The result can be proved by using the dominated-convergence theorem: let
		\begin{align*}
		f_n(\theta) &:= -\tr\log\left(\Psi^{-1}+Q+\frac{1}{n}I_m\right),\quad n\ge 1,\\
		f(\theta) &:= -\tr\log\left(\Psi^{-1}+Q\right),
		\end{align*}
		defined for $\theta\in[-\pi,\pi]$. Clearly $f_n$ is a continuous and therefore measurable function of $\theta$, for any $n\ge1$, and moreover $f_n\uparrow f$ pointwise (recall that $\Psi^{-1}$ and $Q$ are continuous functions of $\theta$) ensuring that  the limit $f$ is itself measurable. In addition, $|f_n|\le g$ pointwise (and therefore a.e.) for any $n\ge1$, where $g=|f|\in L^1(\TT,\HH_m)$ since $Q\in\Q_\Psi^+$. By Lebesgue's dominated-convergence theorem, $f\in L^1(\TT,\HH_m)$ and 
		\begin{align*}
		&\lim_{n\to\infty}\int\,-\tr\log\left(\Psi^{-1}+Q+\frac{1}{n}I_m\right)=\lim_{n\to\infty}\int\,f_n\\
		&= \int\,\lim_{n\to\infty} f_n = \int -\tr\log\left(\Psi^{-1}+Q\right)
		\end{align*}
		where the last equality follows from continuity of $\tr(\cdot)$ and $\log(\cdot)$. The conclusion is now straightforward: 
		\begin{align*}
		\J_\Psi^\infty(Q)&= \lim_{n\to\infty} \J_\Psi^n(Q)=\int\,\tr\left[Q\,\hat{\Phi}_n\right]\\
		&-\lim_{n\to\infty}\int\,\tr\log\left(\Psi^{-1}+Q+\frac{1}{n}I_m\right) +\lambda\,h_\sigma^\infty(Q)\\
		&=\J_\Psi(Q).
		\end{align*}
		for any $Q\in\Q_\Psi^+$. 
\end{IEEEproof}		

We now show our existence and uniqueness result for the extended function $\J_\Psi^\infty$.
\begin{lemma}
The function $\J_\Psi^\infty$ admits a unique minimum on $\bar{\Q}_\Psi^+$.
\end{lemma}
\begin{IEEEproof}
Since we have already seen that  $\J_\Psi^\infty$ is strictly convex on $\bar{\Q}_\Psi^+$, it is sufficient to show that it is also proper on $\bar{\Q}_\Psi^+$. First, notice that $\J_\Psi^\infty$ is not identically $+\infty$ on $\bar{\Q}_\Psi^+$, and it is easy to see that for any $n\ge 1$,
\begin{align}\label{eq:lowbound}
\J_\Psi^n(Q) &\ge \text{cost.} + \int_{-\pi}^{\pi}\, \mu_1\sum_{i=1}^m\,\lambda_i-\frac{1}{\mu_1}\,\log\left(\lambda_i+\frac{1}{n}\right)\,\frac{d\theta}{2\pi}\\
&\ge \text{cost.} +\frac{1}{\mu_1^2}+\frac{1}{n}>-\infty,\notag
\end{align}  
where $\mu_1(e^{i\theta})$ is the minimum eigenvalue of $\hat{\Phi}_n$ while $\lambda_1(e^{i\theta})\le\cdots\le\lambda_m(e^{i\theta})$ are the eigenvalues of $\Psi^{-1}+Q$. Taking the limit for $n\to\infty$ on both sides, we conclude that $\J_\Psi^\infty>-\infty$ on $\bar{\Q}_\Psi^+$. The last step consists in proving that 
\[
\J_\Psi^\infty(Q^{(k)})\to+\infty\qquad\text{when}\qquad(Q^{(k)})_{k\ge1}\subseteq\bar{\Q}_\Psi^+:\,\|Q^{(k)}\|_\mathcal{P}\to\infty.
\]
Let $(Q^{(k)})_{k\ge1}\subseteq\bar{\Q}_\Psi^+$ be such that $\|Q^{(k)}\|_\mathcal{P}\to\infty$ as $k\to\infty$ and denote by $\lambda_1^{(k)}(e^{i\theta})\le\cdots\le\lambda_m^{(k)}(e^{i\theta})$ the eigenvalues of $\Psi^{-1}+Q^{(k)}$ for any $k\ge1$. With $\xi_i(e^{i\theta})$ and $\nu_i^{(k)}(e^{i\theta})$, $i=1,\dots,m$, we denote the eigenvalues of $\Psi^{-1}$ and $Q^{(k)}$ respectively, assumed to be ordered as the $\lambda_i^{(k)}$s. For the sequence $(Q^{(k)})_{k\ge1}$ we have $\|Q^{(k)}\|_\mathcal{P}\to\infty$ as $k\to\infty$. Starting from \eqref{eq:lowbound} and applying Weyl's Theorem \cite[Section 4.3]{horn1990} we get
\begin{equation}\label{eq:ineqJn1}
\scalebox{0.9}{$\displaystyle
\begin{split}
\J_\Psi^n(Q^{(k)})&\ge \text{cost.} + \int\, \left[\mu_1\sum_{i=1}^m\,(\xi_1+\nu_i^{(k)})-\frac{1}{\mu_1}\,\log\left(\nu_i^{(k)}+\xi_m+\frac{1}{n}\right)\right]\\
&\ge \text{cost.} + \mu_1\,\int\,\left[ \xi_m+\nu_m^{(k)}-\frac{m}{\mu_1}\,\log\left(\nu_m^{(k)}+\xi_m+\frac{1}{n}\right)\right]\\
&\ge \text{cost.} + \mu_1\,\int\,\left[ \nu_m^{(k)}\right]-\frac{m}{\mu_1}\,\log\left(\int\, \nu_m^{(k)}+\xi_m+\frac{1}{n}\right),
\end{split}
$} 
\end{equation}
where the last step follows from the Jensen's inequality. Now observe that
\begin{equation}\label{eq:normdiv}
\|Q^{(k)}\|_\mathcal{P}=\int\,|\nu_m^{(k)}|=\int_{\V_+^{(k)}}\nu_m^{(k)}+\int_{\V_-^{(k)}}(-\nu_m^{(k)}),
\end{equation}
where
\begin{align*}
\V_+^{(k)}&:=\left\{\theta\in[-\pi,\pi]:\,\nu_m^{(k)}(e^{i\theta})>0\right\},\\
\V_-^{(k)}&:=\left\{\theta\in[-\pi,\pi]:\,\nu_m^{(k)}(e^{i\theta})<0\right\}.
\end{align*}
From Weyl's Theorem $-\nu_m^{(k)}\le\xi_m$, therefore the second integral in \eqref{eq:normdiv} is bounded above and
\[
\|Q^{(k)}\|_\mathcal{P}\to\infty\qquad\implies\qquad \int_{\V_+^{(k)}}\nu_m^{(k)}\to\infty.
\]
From \eqref{eq:ineqJn1} the inequality
\begin{align*}
\J_\Psi^n&(Q^{(k)}) \ge \text{cost.} + \mu_1\,\int\, \left[\nu_m^{(k)}\right]-\mu_1\,\int_{\V_-^{(k)}}\left[-\nu_m^{(k)}\right]\\
&-\frac{m}{\mu_1}\,\log\left(\int\left[\nu_m^{(k)}\right]-\int_{\V_-^{(k)}}\left[-\nu_m^{(k)})\right]+\xi_m+\frac{1}{n}\right)
\end{align*}
holds for any $n\ge 1$. Taking the limit for $k\to+\infty$ on both sides we obtain that $\J_\Psi^n(Q^{(k)})\to+\infty$ when $\|Q^{(k)}\|_\mathcal{P}\to\infty$ for any $n\ge1$ and therefore also when $n\to\infty$, i.e. $\J_\Psi^\infty(Q^{(k)})\to+\infty$ when $\|Q^{(k)}\|_\mathcal{P}\to\infty$. Given the fact that the index $\J_\Psi^\infty$ is strictly convex over $\bar{\Q}_\Psi^+$, a Weierstrass'-theorem argument \cite[pp. 35]{ekeland1999} allows to conclude that $\J_\Psi^\infty$ admits \emph{unique} minimum point $\bar{\Q}_\Psi^+$.
\end{IEEEproof} 

The last step consists in showing that the minimum of $\J_\Psi^\infty$ cannot be on the boundary.
\begin{lemma}\label{lemma-interior}
The minimum of $\J_\Psi^\infty$ is attained in $\Q_\Psi^+$.
\end{lemma}

\begin{IEEEproof}
Let $Q_o\in\partial\Q_\Psi^+$. For any $\epsilon>0$, $Q_o+\epsilon I_m\in\bar{\Q}_\Psi^+$ and  $h_\sigma^\infty(Q_o+\epsilon I_m)=h_\sigma^\infty(Q_o)$. Bringing the limit inside the integral, we can find an upper-bound to the right G\^ateaux derivative of $\J_\Psi^\infty$ in direction $\delta Q=I_m$, i.e.
		\[
		\delta \J_\Psi^\infty(Q_o;I_m)=\lim_{\epsilon\downarrow 0} \frac{\J_\Psi^\infty(Q_o+\epsilon I_m)-\J_\Psi^\infty(Q_o)}{\epsilon}.
		\]
		Notice that $\delta \J_\Psi^\infty(Q_o;I_m)=dF/d\epsilon$ where 
		\[
		F(\epsilon):=\int_{-\pi}^\pi\,\tr\left[\left(Q_o+\epsilon I\right)\hat{\Phi}_n-\log\left(\Psi^{-1}+Q_o+\epsilon I\right)\right]\,\frac{d\theta}{2\pi},
		\]
		defined for any $\epsilon>0$. In fact, notice that 
		\[
		f(\theta,\epsilon):=\tr\left[\left(Q_o+\epsilon I\right)\hat{\Phi}_n-\log\left(\Psi^{-1}+Q_o+\epsilon I\right)\right]
		\]	
		is integrable for each $\epsilon>0$ and, by a similar argument as the one used in \eqref{eq:ineqJn1},
		\[
		\bigg|\frac{\partial f(\theta,\epsilon)}{\partial\epsilon}\bigg|\le \big|\tr\left(\hat{\Phi}_n\right)\big|+(\lambda_m+1)\big|\tr\left(\Psi^{-1}+Q_o+\epsilon I\right)\big|=:g(\theta),
		\]
		with $g\in L^1(\TT,\HH_m)$, for any $\epsilon>0$ and $\theta\in[-\pi,\pi]$. Accordingly, $F$ is differentiable and
		\begin{align*}
		\delta\J_\Psi^\infty(Q_o;I_m)&=\frac{dF(\epsilon)}{d\epsilon}= \int\frac{\partial f(\theta,\epsilon)}{\partial\epsilon}\\
		&=\int\tr\left(\hat{\Phi}_n\right)-\tr\left[I+\left(\Psi^{-1}+Q_o\right)^{-1}\right]\\
		&\le m\int\,\mu_m\,\frac{d\theta}{2\pi} -\int \tr\left[I+\left(\Psi^{-1}+Q_o\right)^{-1}\right]\\
		\end{align*}
	where $\mu_m(e^{i\theta}),\,\theta\in[-\pi,\pi]$, is the maximum eigenvalue of $\hat{\Phi}_n$. Inasmuch $Q_o\in\partial\Q_\Psi^+$, $\tr(\Psi^{-1}+Q_o)^{-1}$ is a positive rational functions having poles on $\TT$, accordingly $\int\,\tr(\Psi^{-1}+Q_o)^{-1}\to+\infty$ and therefore $\delta \J_\Psi^\infty(Q_o;I_m)\to-\infty$. Hence, for $\epsilon>0$, sufficiently small, $\J_\Psi^\infty(Q_o+\epsilon I)<\J_\Psi^\infty(Q_o)$ so that $Q_o$ cannot be the minimum. Equivalently, the minimum cannot belong to $\partial\Q_\Psi^+$.
\end{IEEEproof}
As a result of the previous steps, we have the following result that provides a solid theoretical ground to our work:
\begin{theorem}\label{thm:exuqDual}
Problem \eqref{eq:regISminDual} admits a unique solution $Q_o\in\Q_\Psi^+$.
\end{theorem}

\medskip

Before concluding the section, it is worth noticing that Problem \eqref{eq:regISminDual} is a regularized (and thus more complex) version of Problem \eqref{eq:ISminDual}. Thus, by following the same steps 
of the proof of Theorem \ref{thm:exuqDual}, we can extend also the objective function of Problem \eqref{eq:ISminDual} to the closure of the original domain $\Q_\Psi^+$. Hence, we can apply Weierstrass Theorem and conclude that Problem \eqref{eq:ISminDual} admits a minimum in the closed domain. Then, by resorting to the same argument of Lemma \ref{lemma-interior} we can show  that the minimum is actually attained on $\Q_\Psi^+$ i.e. it cannot be achieved on the boundary.
Finally, by duality occurring between Problem \eqref{eq:ISmin} and Problem \eqref{eq:ISminDual}, 
we have the following ancillary  result that, as discussed in Remark \ref{remark-avventi},
solves the problem  of identifying  an ARMA graphical model with topology $\Omega_\tau$ thus extending to a general prior the work \cite{avventi2013}.

\begin{theorem}\label{thm:exuqDemp}
Problem \eqref{eq:ISmin} admits a unique solution. 
\end{theorem}

\section{Simulation Results}\label{sec:sim}
In this section we present some numerical examples illustrating the performances of the proposed method for positive link prediction problems. More specifically, the experiments test the proposed algorithm over a network of agents in two directions: the case in which the underlying dynamics can be well approximated by an AR model, i.e. the involved spectra are trigonometric polynomials, and the case in which the approximation is made through an ARMA model.

\subsubsection*{AR dynamics} 
\begin{figure}[h!]
	\centering
	\includegraphics[scale=0.4]{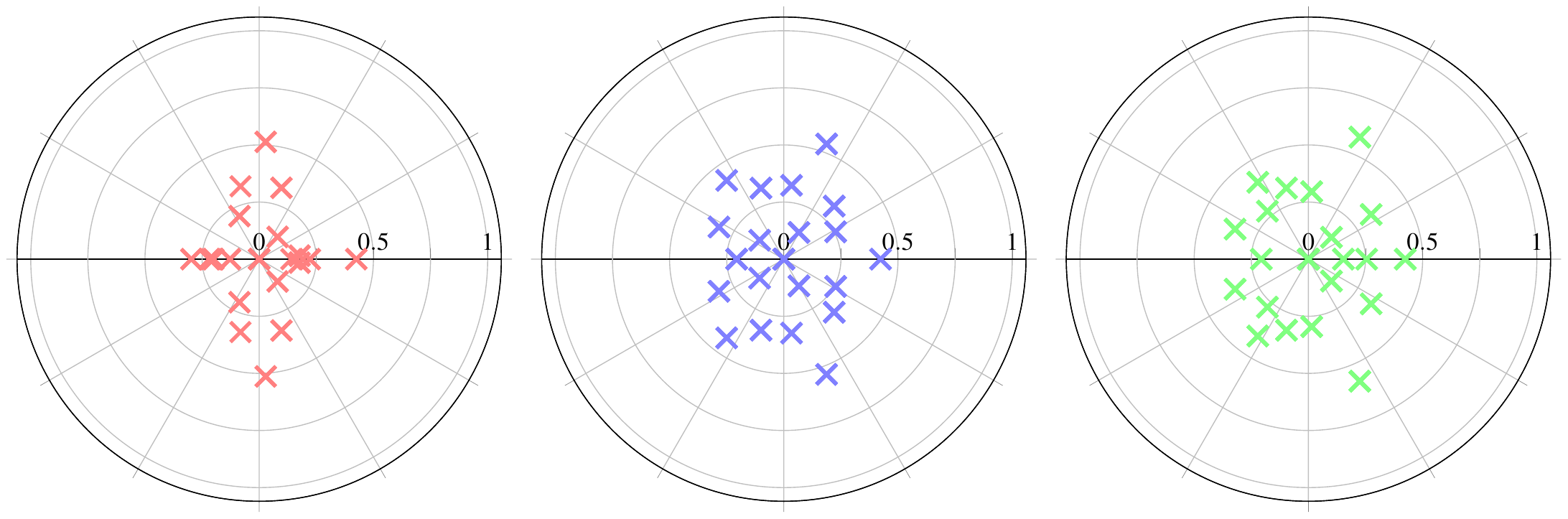}
	\caption{Poles the shaping filter $[I-F(z)]^{-1}$ of the AR model: time $\tau=0$ (left), time $\tau=1$ (center), time $\tau=2$ (right).}
	\label{fig:poles_AR}
\end{figure}

\begin{figure}[h!]
	\centering
	\includegraphics[scale=0.7]{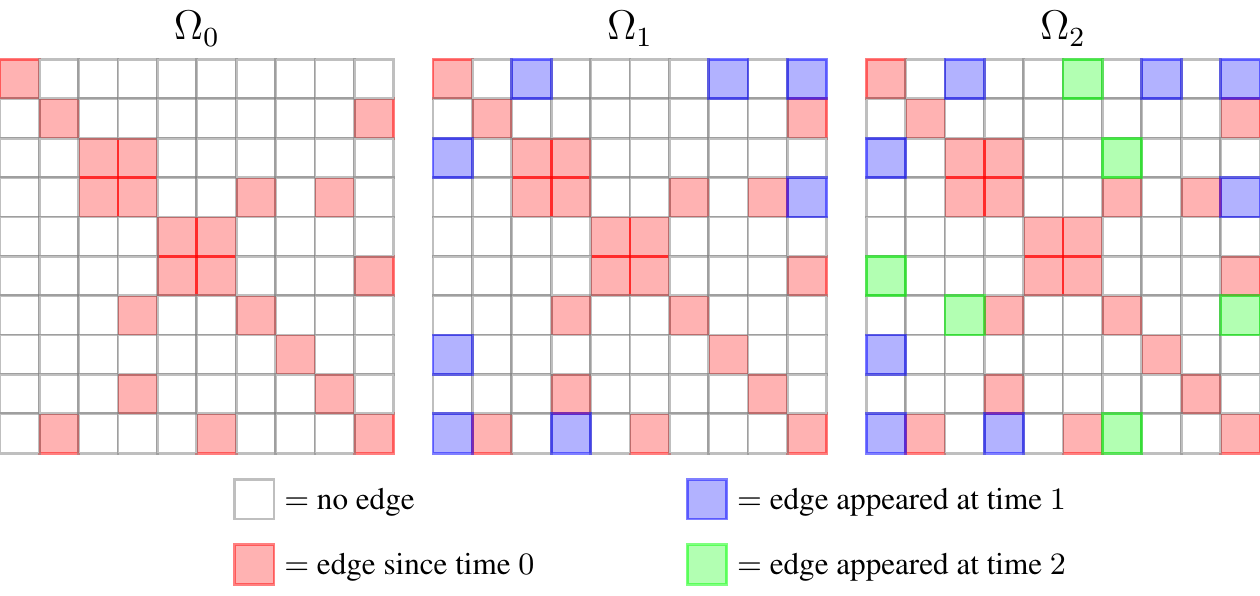}
	\caption{Support of the prior $\Phi_0^{-1}$ (left), and true supports of the spectra $\Phi_1^{-1}$ (center) and $\Phi_2^{-1}$ (right) that have to be estimated.}
	\label{fig:true_supp_PLP_AR}
\end{figure}

Suppose that the dynamic of the agents composing the network is described by model \eqref{eq:ARMA} with $B_k=0$ and $F_k=-A_0^{-1}A_k$ for all $k>0$,
\begin{equation}\label{eq:ARsim}
    \rv{y}(t) = \sum_{k=1}^{n}\,F_k\,\rv{y}(t-k) + \rv{e}(t),
\end{equation}
i.e. by the AR process $\rv{y}(t) = F(z)\,\rv{y}(t) + \rv{e}(t)$ in which $F(z)=\sum_{k=0}^nF_k\,z^{-k}$, the dimension of the process is $m=10$ and the order of the process is $n=2$. The set-up for this test is the one of Recursive PLP explained in Section \ref{sec:recursive}, for a window of length $T=2$. The initial information of the network is enclosed in a prior spectral density $\Phi_0$ whose inverse has support $\Omega_0$ and the supports $\Omega_1$ and $\Omega_2$ of the spectra $\Phi_1^{-1}$ and $\Phi_2^{-1}$ need to be estimated according to the iterative scheme
\[
   \Phi_0 \leadsto \Phi_1 \leadsto \Phi_2,
\]
outlined in Section \ref{sec:recursive}. Figure \ref{fig:poles_AR} depicts the poles of the shaping filter $[I-F(z)]^{-1}$ of model \eqref{eq:ARsim} at times $\tau=0,1,2$ from left to right respectively, while Figure \ref{fig:true_supp_PLP_AR} reports the support of $\Phi_0^{-1}$ (left) and the supports of the spectra $\Phi_1^{-1}$ (center) and $\Phi_2^{-1}$ (right) that have to be estimated. It is worth noting that we are dealing with an unfriendly prediction network. Indeed, being $\N_i$ the set of neighbors of agent $i$ according to the network topology of $\G(0)$, consider the common neighbors similarity measure at time $\tau=1$,
\[
   (CN_0)_{ij}=\mathrm{card}(\N_i\cap \N_j),
\]
where $\mathrm{card}(\N_i\cap \N_j)$ denotes the cardinality of set $\N_i\cap \N_j$. It is not difficult to see that $(CN_0)_{ij}=0$ for any edge $(i,j)$ appearing in $\G(1)$, while $(CN_0)_{ij}=1$ for any $(i,j)\in\Omega^{CN}_1:=\{(2,6),\, (3,7),\,(3,9)\,(5,10),(7,9) \}$ and $\Omega^{CN}_1$ contains edges not appearing in $\G(1)$. A similar scenario happens for the prediction time $\tau=2$ revealing that the common neighbors measure leads to low prediction accuracy for this kind of network.
\\%

\begin{figure}[h!]
	\centering
	\includegraphics[scale=1.65]{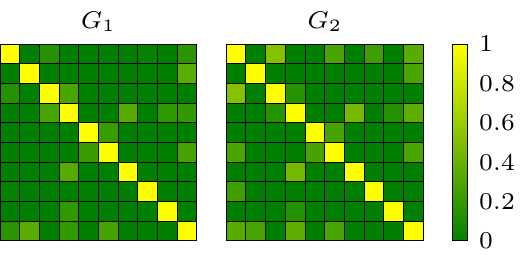}
	\caption{Score matrices $G_1$ and $G_2$ for estimation of the supports $\Omega_1$ and $\Omega_2$, respectively.}
	\label{fig:scores}
\end{figure}

\begin{figure}[h!]
	\centering
	\includegraphics[scale=0.65]{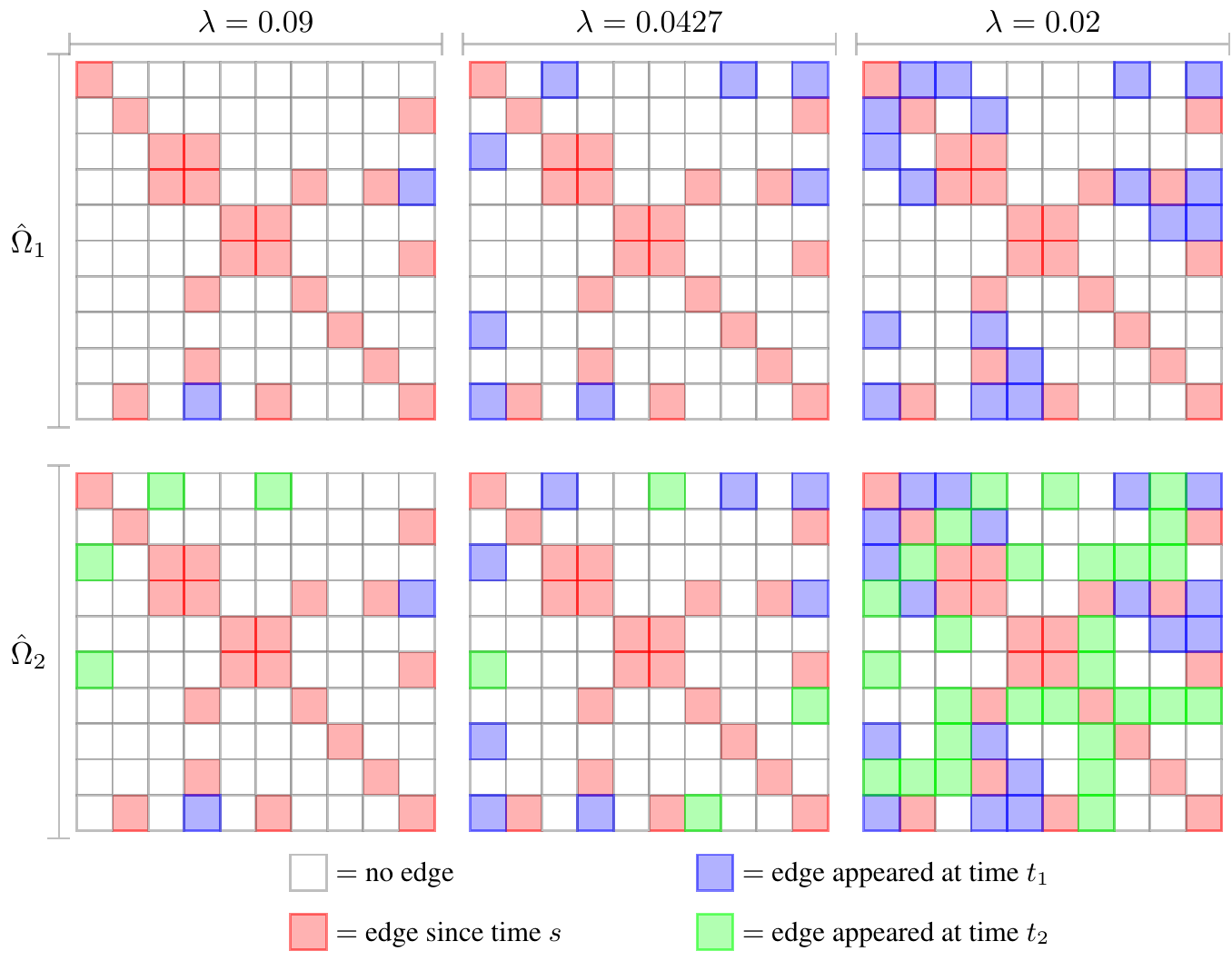}
	\caption{Estimates of the supports $\Omega_1$ and $\Omega_2$ for different values of the regularization parameter $\lambda$.}
	\label{fig:est_supp_PLP_AR}
\end{figure}

In regard to our method, the computation of the estimates $\hat{\Omega}_1$ and $\hat{\Omega}_2$ of the respective supports $\Omega_1$ and $\Omega_2$, exploits the information coming from $N_1=N_2=1000$ data samples together with the previously estimated spectrum and it is performed by solving Problem \eqref{eq:recopt} for each estimation step using the CVX package of Matlab \cite{GrantBoyd08,GrantBoyd14}. As mentioned in Section \ref{sec:PLP} the estimates of the inverse spectra produced by our method are close to be sparse. Figure \ref{fig:scores} displays the score matrices $G_\tau$, $\tau=1,2$, obtained by evaluating the partial coherence-based similarity measure with regularization parameter $\lambda=0.0427$ that has been kept constant for both the estimation steps.
\\%
Accordingly, the estimates of the supports, namely the estimated network topologies, have been obtained by a thresholding procedure that sets to zero all the entries in position $(i,j)\in V\times V\setminus\Omega_\sigma$ such that $(G_\tau)_{ij}<t_r$. The threshold value for this simulation is $t_r=3\cdot10^{-4}$, kept constant for both the estimation steps.\\
Figure \ref{fig:est_supp_PLP_AR} displays the resulting support estimates for different values of the regularization parameter $\lambda$, that has been kept constant for both the estimation steps. The value $\lambda=0.0427$ is a good choice in that the procedure has reached the perfect recovery of $\Omega_1$ and only one non-zero entry is missing in the estimate $\hat{\Omega}_2$ when compared to $\Omega_2$. To give a wider view on the performances of the proposed algorithm, Figure \ref{fig:est_supp_PLP_AR} reports also the results for $\lambda=0.09$ (left column), which results in excessively sparse estimates, and the results for $\lambda=0.02$ (right column) in which a weak regularization effect is highlighted.

\subsubsection*{ARMA dynamics} 
\begin{figure}[h!]
	\centering
	\includegraphics[scale=0.7]{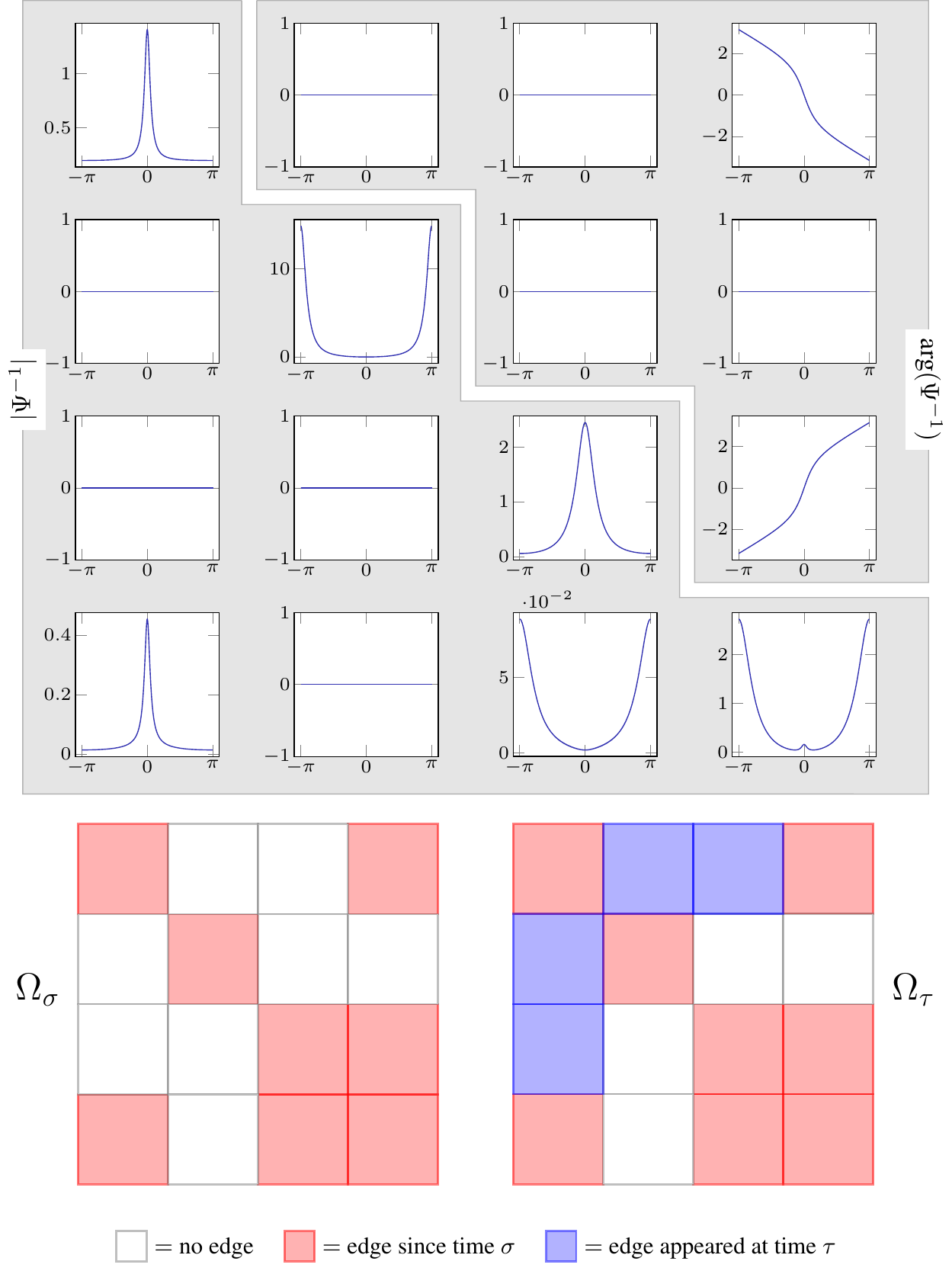}
	\caption{Inverse of the ARMA prior (above), prior's support (below left) and true support of $\Phi^{-1}$ (below right).}
	\label{fig:true_supp_PLP_ARMA}
\end{figure}

\begin{figure}[h!]
	\centering
	\includegraphics[scale=0.7]{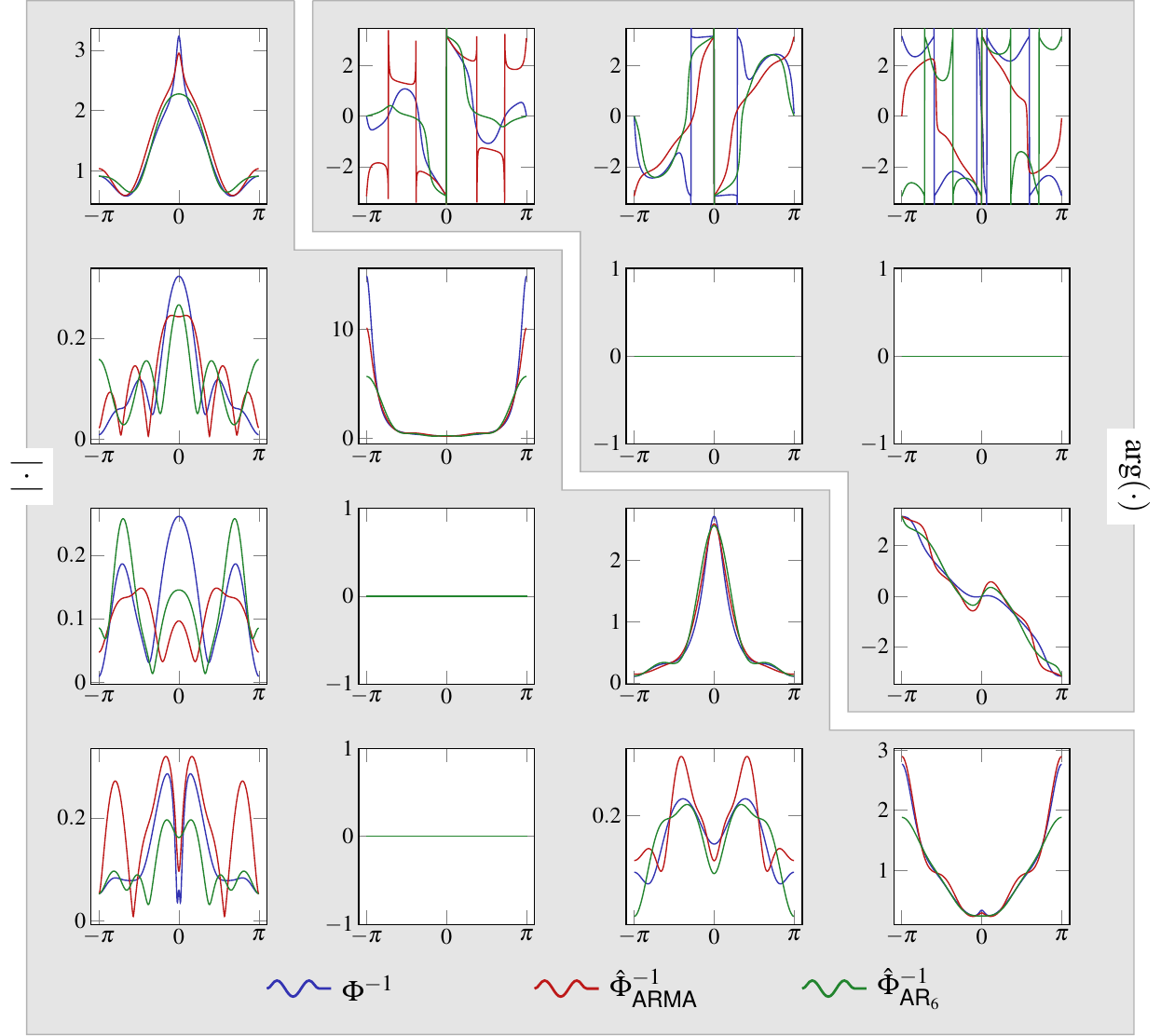}
	\caption{Comparison between the true $\Phi^{-1}$ and its estimate $\hat{\Phi}_{\textsf{ARMA}}^{-1}$ computed from the ARMA prior and its estimate $\hat{\Phi}_{\textsf{AR}_6}^{-1}$ computed with no prior information.}
	\label{fig:Psds_plot}
\end{figure}

This test illustrates the performances of the proposed method when the model used to fit the data is an ARMA model of the type of \eqref{eq:ARMA} in which the dimension of the process is $m=4$ while the order of the polynomial part $Q$ is set to $n=4$. This experiment generalizes the approach proposed in \cite{SongVan2010} to the case of a vector MA part. To relate our results with \cite{SongVan2010}, we actually consider the solution of Problem \eqref{eq:regISminDual} in two different cases: the case in which an ARMA prior spectral density $\Psi^{-1}$ and its support $\Omega_\sigma$ are available, and the case in which $\Psi\equiv I_m$ (and $\Omega_\sigma=I_m$) namely, the regularized maximum likelihood estimator of \cite{SongVan2010}, in which we require our estimation procedure to fit the ARMA model \eqref{eq:ARMA} with an AR model of order $n=6$ (so that to have a comparable number of parameters). Figure \ref{fig:true_supp_PLP_ARMA} reports the inverse of the ARMA prior power spectral density chosen for the first test (above) with its support $\Omega_\sigma$ (below left), and the support $\Omega_\tau$ of $\Phi^{-1}$ that has to be estimated (below right). 
The results of the simulations corresponding to $\lambda=0.04$ are illustrated in Figure \ref{fig:Psds_plot} where the true $\Phi^{-1}$ (blue) is compared to its estimates obtained from the two different priors. In particular, the red line represents the estimate $\hat{\Phi}_{\textsf{ARMA}}^{-1}$ of $\Phi^{-1}$ computed from the ARMA prior $\Psi$ depicted in Figure \ref{fig:true_supp_PLP_ARMA} while the estimate $\hat{\Phi}_{\textsf{AR}_6}^{-1}$ computed without prior information is represented with the green line. 
\\%
Also for these set of simulations, the supports of the estimated spectra has been computed on the basis of the partial coherence-based similarity measure followed by the thresholding procedure (with $t_r=0.1$) and we can see from Figure \ref{fig:Psds_plot} how the estimation procedure is able to recover the true support $\Omega_\tau$ with both priors. However, it can be noticed that the presence of the prior leads to an enhancement of the estimation capabilities of the proposed procedure. This is particularly highlighted by the comparison between the entries in position $(4,1)$ in Figure \ref{fig:Psds_plot}. The proposed algorithm has been tested also in the case in which the true model has been tested also on an AR model of the same order $n=4$ of the ARMA model used for the estimation. As expected, the resulting estimate $\hat{\Phi}_{\textsf{AR}_4}$ has turned out to be worse than the estimate $\hat{\Phi}_{\textsf{AR}_6}$.

\section{CONCLUSIONS}\label{sec:conc}
In this work, positive link prediction as detection problem is approached as an identification problem for ARMA graphical models when some a-priori information is available. The main contribution of this work is the introduction of a similarity measure that instead of relying on  properties that the network is expected to fulfill, it relies on noisy observations of the network at the current time, setting the proposed link-detection method in the context of (partially) data-driven approaches. The positive link prediction problem was rephrased as a suitable optimization problem whose solution has been formally proved to exist and to be unique. Although this work is mainly theoretically-focused, the solution has been computed numerically for different synthetic-data case studies and the method has been compared with the existing methods for the identification of graphical models showing an improvement of the performances with respect to the case in which no a-priori information is available.


\end{document}